\documentclass{ifacconf}

\usepackage{graphicx}      
\usepackage{xcolor}

\makeatletter
\let\old@ssect\@ssect 
\makeatother
\usepackage{natbib}
\usepackage{hyperref}
\hypersetup{
   breaklinks=true,   
   pdfusetitle=true,  
   hidelinks
}
\makeatletter
\def\@ssect#1#2#3#4#5#6{%
  \NR@gettitle{#6}
  \old@ssect{#1}{#2}{#3}{#4}{#5}{#6}
}
\makeatother

\usepackage{mathtools}
\usepackage{amsmath}
\usepackage{amssymb}
\usepackage{enumitem}
\usepackage{booktabs} 

\newcommand{\Dim}[1]{\mathrm{dim}(#1)}

\newcommand{\D}{\mathrm{d}}

\newcommand{\X}{\mathcal{X}}
\newcommand{\Qm}{\mathcal{Q}}
\newcommand{\T}{\mathcal{T}}

\newcommand\undermat[2]{%
	\makebox[0pt][l]{$\smash{\underbrace{\phantom{%
					\begin{matrix}#2\end{matrix}}}_{\text{$#1$}}}$}#2}
\begin{document}
\begin{frontmatter}

\title{ Exact Linearization of Minimally Underactuated Configuration Flat Lagrangian Control Systems by Quasi-Static Feedback of Classical States } 

\thanks[footnoteinfo]{ This research was funded in whole, or in part, by the Austrian Science Fund (FWF) P36473. For the purpose of open access, the author has applied a CC BY public copyright licence to any Author Accepted Manuscript version arising from this submission. }

\author[First]{Georg Hartl} 
\author[First]{Conrad Gstöttner} 
\author[First]{Bernd Kolar}
\author[First]{Markus Schöberl}

\address[First]{Institute of Automatic Control and Control Systems Technology, Johannes Kepler University Linz, Altenberger Strasse 69, 4040 Linz, Austria (e-mail: \{georg.hartl, conrad.gstoettner, bernd.kolar, markus.schoeberl\}@jku.at).}

\begin{abstract}                

We study the exact linearization of configuration flat Lagrangian control systems with $p$ degrees of freedom and $p-1$ inputs by quasi-static feedback of classical states. First, we present a detailed analysis of the structure of the parameterization of the system variables by the flat output. Based on that, we systematically construct a linearizing quasi-static feedback law of the classical state such that the closed-loop system shows the behavior of decoupled integrator chains. Our approach shows that the construction of a generalized Brunovský state can be completely circumvented. Furthermore, we present a method for determining the lengths of the integrator chains achieved by quasi-static feedback laws that allow for rest-to-rest transitions.

\end{abstract}


\end{frontmatter}

\section{ Introduction }
One of the most significant contributions to nonlinear control theory in recent decades has been the introduction of the concept of differential flatness by Fliess, Lévine, Martin, and Rouchon, see, e.g., \cite{FliessLevineMartinRouchon:1995}. Differential flatness is a fundamental property of dynamical systems, allowing for systematic solutions in control engineering problems such as trajectory planning or feedback control. Considering nonlinear control systems of the form 
\begin{equation}\label{sys_non_lin_def}
	\dot{x} = f(x,u),
\end{equation}
with $\Dim{x}=n$ states and $\Dim{u}=m$ control inputs, flatness is defined by the existence of an $m$-tuple $y$ of differentially independent functions $y_j=\varphi_j(x,u,u^{(1)}, \ldots, u^{(\nu)})$, where $u^{(\nu)}$ denotes the $\nu$-th time derivative of $u$, that allow expressing the state $x$ and the input $u$ in terms of the so-called flat output $y$ and a finite number of its time derivatives. Note, however, that the computation of flat outputs remains a complex challenge with no comprehensive solution as discussed in recent research, including, e.g., \cite{Gstottner:2023-2}, \cite{NicolauRespondek:2017}, and \cite{SchoberlSchlacher:2014}. The underlying principle of flatness-based control is the so-called exact linearization, where a nonlinear system (\ref{sys_non_lin_def}) is exactly linearized by endogenous feedback such that the closed-loop system has a linear input-output behavior 
\begin{equation}
	\label{in_out}
	y_j^{(\kappa_j)} = w_j, \quad j=1,\ldots,m,
\end{equation}
between a new input $w$ and the flat output $y$, in form of $m$ decoupled integrator chains with respective lengths $\kappa_j$. 

As shown, e.g., in \cite{FliessLevineMartinRouchon:1999}, every flat system is linearizable by a systematically derivable endogenous dynamic feedback law. However, applying such feedback laws increases the order of the tracking error dynamics. A widely acknowledged alternative introduced, e.g., in \cite{DelaleauRudolph:1995} and \cite{DelaleauRudolph:1998}, is the exact linearization by quasi-static feedback of so-called generalized Brunovský states \mbox{$\Tilde{x}_B=(y, \dot{y}, \ldots, y^{(\kappa - 1)})$}. Every flat system possesses such generalized Brunovský states, which consist of suitably chosen consecutive time derivatives of the flat output $y$ and meet $\Dim{\Tilde{x}_B}=\Dim{x}$. Thus, a flat system (\ref{sys_non_lin_def}) can be exactly linearized by a quasi-static feedback law \mbox{$u = \alpha(\Tilde{x}_B,w, w^{(1)}, \ldots)$}, where the new input $w$ is given by suitable time derivatives $y^{(\kappa)}$ of the flat output. 

Contrary to dynamic feedback, quasi-static feedback does not increase the order of the closed-loop system dynamics but may include time derivatives of the new input $w$ of the closed-loop system. However, calculating the required time derivatives of the flat output based on the available measurements to determine $\Tilde{x}_B$ is generally challenging. To address the issue of determining $\Tilde{x}_B$, e.g., \cite{KolarRamsSchlacher:2017} demonstrated that, under certain conditions, $\Tilde{x}_B$ can be replaced by the classical state $x$ in the quasi-static feedback. Using a differential algebraic setting, \cite{DelaleauFiltrationsFeedbackSynthesis1998} demonstrated that every system with a flat output that depends only on the classical state $x$ and the input $u$, so-called $(x,u)$-flat systems, can be exactly linearized by a quasi-static feedback law $u=\alpha(x,w,w^{(1)},\ldots)$ of the classical state $x$. Utilizing a differential geometric framework, similar results for the exact linearization of $(x,u)$-flat systems are presented by \cite{Gstottner:2023-1}, where a systematic procedure for constructing a quasi-static feedback law of the classical state $x$ is presented. In \cite{Gstottner:2023-1}, the multi-index $\kappa$, which contains the orders of the time derivatives of the flat output $y$ intended as new inputs, is constructed utilizing an algorithm. Through successive coordinate transformations within a finite-dimensional differential geometric framework, this algorithm provides the multi-index $\kappa$ and the quasi-static feedback of the classical state, all achieved without the use of generalized Brunovský states. For flat discrete-time systems, an analogous approach to \cite{Gstottner:2023-1} is presented in \cite{KolarDiwoldGstottnerSchoberl:2022}.



It has been shown that many mechanical systems are actually flat. An analysis of the flatness of mechanical systems with three degrees of freedom is presented, e.g., in \cite{Nicolau:2015}. In the present work, we restrict ourselves to Lagrangian control systems with a Lagrangian of the form kinetic minus potential energy as, e.g., in \cite{RathinamMurray:1998} and \cite{SatoIwai:2012}. Specifically, \cite{RathinamMurray:1998} present necessary and sufficient conditions for Lagrangian control systems with $p$ degrees of freedom and $p-1$ control inputs to be configuration flat, i.e., possessing a flat output that depends only on configuration variables.   

This paper is devoted to the exact linearization of configuration flat Lagrangian control systems with $p$ degrees of freedom and $p-1$ control inputs, i.e., minimally underactuated systems, by quasi-static feedback of classical states. Thereby, we circumvent entirely the use of generalized Brunovský states. Note that Lagrangian control systems with $p$ degrees of freedom and $p$ inputs, i.e., fully actuated systems, are generally flat and always exactly linearizable by static state feedback. Systems with $p-1$ inputs represent the next more complicated case, covering various practical examples, such as cranes or aircraft. In the context of flatness, such systems have been examined, e.g., in \cite{RathinamMurray:1998} or \cite{FliessLevineMartinRouchon:1999}. 

The main idea of our approach is based on essential findings of \cite{Gstottner:2023-1}. However, we do not utilize successive coordinate transformations in a differential geometric framework. Instead, we present a detailed analysis of the structure of the parameterization of the system variables by the flat output $y$. Subsequently, we demonstrate that all possible orders $\kappa_j$ of the time derivatives of the flat output $y$ intended as new inputs are easily identifiable by the Jacobian matrix of the parameterization of the classical state by the flat output $y$. Typically, control engineering tasks for mechanical systems involve trajectory tracking for transitions from one equilibrium point to another, so-called rest-to-rest transitions. We show that linearizing feedback laws that inevitably have singularities at equilibrium points can be formerly excluded. Finally, we derive a quasi-static feedback law of the classical state. This work is structured as follows: Section \ref{sec_2} introduces the notation used throughout this paper. In Section \ref{sec_4}, we present properties of differentially flat systems and introduce configuration flat Lagrangian control systems. Section \ref{sec_5} is dedicated to our main results, illustrated by one practical example in Section \ref{sec_6}.

\section{Notation}
\label{sec_2}
We use tensor notation and the Einstein summation convention if the index range is clear from the context. Consider an $n$-dimensional smooth manifold $\mathcal{M}$ with $x$ denoting local coordinates $(x^1,\dots,x^n)$ and an $m$-tuple of functions on $\mathcal{M}$, $f=(f^1,\dots,f^m): \mathcal{M} \rightarrow \mathbb{R}^m$. We denote the $m \times n$ Jacobian matrix by $\partial_x f$ and the partial derivative of the function $f^j$ with respect to the variable $x^i$ by $\partial_{x^i} f^j$. Likewise, $\text{d}f$ stands for the differentials $(\text{d}f^1,\dots,\text{d}f^m)$. Moreover, we express time derivatives of variables, for example, the control inputs $u$, and functions, for example, the flat output $\varphi$, by subscripts in squared brackets instead of employing the notation with dots and superscripts in round brackets as used, e.g., in \cite{FliessLevineMartinRouchon:1999}. For instance, $y^j_{[\alpha]}$ stands for the $\alpha$-th time derivative of the $j$-th component of the $m$-tuple $y$. Further, $y_{[\alpha]} = (y^1_{[\alpha]}, \dots , y^m_{[\alpha]})$ denotes the $\alpha$-th time derivative of every component of $y$. Multi-indices marked by capital letters denote time derivatives of different orders of each component of a tuple $y=(y^1, \dots, y^m)$. Let $A = (a^1, \dots , a^m)$ and $B=(b^1, \dots ,b^m)$ be two multi-indices with $a^j \leq b^j$ for all $j=1, \dots, m$, which is abbreviated by $A \leq B$. Likewise, $A\leq c$, where $c$ is an integer, signifies that $a^j \leq c$ for all $j=1, \dots, m$. By using these multi-indices, we introduce the abbreviations \mbox{$y_{[A]} = (y^1_{[a^1]}, \dots ,y^m_{[a^m]})$}, \mbox{$y_{[0,A]} = ( y^1_{[0,a^1]}, \dots ,y^m_{[0,a^m]} )$} and \mbox{$y_{[A,B]} = ( y^1_{[a^1,b^1]}, \dots ,y^m_{[a^m,b^m]} )$}, where the entries of the last tuple represent \mbox{$y^j_{[a^j, b^j ]} = ( y^j_{[a^j]}, \dots , y^j_{[b^j]})$}. The addition and subtraction of multi-indices is done component-wise by $A \pm B = ( a^1\pm b^1, \dots, a^m \pm b^m )$. We denote the summation over the indices by \mbox{$\#A = \sum_{j=1}^m a^j$}.

\section{ Preliminaries }\label{sec_4}
\vspace{-3pt}
\subsection{ Differential flatness and exact linearization }
\vspace{-3pt}
This section provides a short review of flat nonlinear control systems (\ref{sys_non_lin_def}) using a finite-dimensional differential geometric framework as, e.g., in \cite{KolarSchoberlSchlacher:2016-3-ifac-online}. We introduce a so-called state-input manifold \mbox{$\X \times \mathcal{U}_{[0,l_u]}$} with local coordinates $(x,u,u_{[1]}, \ldots, u_{[l_u]})$ and some large enough order $l_u$ to define flatness as follows:
\begin{defn}\label{def_1}
    A system (\ref{sys_non_lin_def}) is called flat if there exists an $m$-tuple of smooth functions $y^j=\varphi^j(x,u,u_{[1]},\ldots), j=1,\ldots,m$, defined on $\X \times \mathcal{U}_{[0,l_u]}$ and smooth functions $F^i_x$ and $F^j_u$ such that locally
    \begin{equation}
        \label{para_map_def}
        \begin{aligned}
            x^i & = F^i_x(y_{[0,R-1]}), \quad && i=1,\ldots, n \\
            u^j & = F^j_u(y_{[0,R]}),  && j=1,\ldots, m, 
        \end{aligned}
    \end{equation}
    with some multi-index $R=(r^1,\ldots,r^m)$. The $m$-tuple $y$ is called a flat output of (\ref{sys_non_lin_def}), and  (\ref{para_map_def}) is called the parameterizing map.
\end{defn}
Throughout this work, we assume that the flat parameterization (\ref{para_map_def}) is well defined locally around equilibrium points. These equilibrium points are characterized by constant values $(x_s, u_s)$, which are in turn parameterized by constant values of the flat output denoted by $y_s$ with $y_{s,[1]}=y_{s,[2]}=\ldots=0$. One implication of Definition \ref{def_1} is that the differentials $\D \varphi, \D \varphi_{[1]}, \ldots, \D \varphi_{[\beta]} $ of derivatives of the flat output up to an arbitrary order $\beta$ are linearly independent. Another implication is that the parameterizing map (\ref{para_map_def}) is unique with a minimal multi-index $R=(r^1, \ldots, r^m)$, where $r^j$ stands for the highest time derivative of $y^j$ required for expressing $x$ and $u$ in terms of the flat output $y$ and its time derivatives. We call systems $(x,u)$-flat if they possess a flat output that depends only on classical states $x$ and inputs $u$ but no time derivatives of $u$. Note that $(x,u)$-flat outputs do not necessarily have to depend explicitly on $u$. Thus, configuration flat systems also belong to the class of $(x,u)$-flat systems.

An intermediate step towards flatness-based tracking control is the exact linearization of flat systems. To exactly linearize a nonlinear system, it is necessary to construct an endogenous feedback law such that the closed-loop system possesses a linear input-output behavior of the form 
\begin{equation}
	\label{input_output}
	y^j_{[\kappa^j]} = w^j, \quad j=1,\ldots,m,
\end{equation}
with a new input $w$ consisting of time derivatives of suitable orders $\kappa=(\kappa^1, \ldots, \kappa^{m})$ of the flat output $y$. When designing flatness-based tracking control, it is advantageous to choose the orders $\kappa^j$ as low as possible to minimize the order of the tracking error dynamics, which can be accomplished by a linearizing quasi-static feedback law. \cite{Gstottner:2023-1} provide easily verifiable conditions for selecting time derivatives of the flat output such that they can be introduced as a new input $w$ by an endogenous feedback law of the classical state. An abbreviated version of these conditions regarding quasi-static feedback of the classical state is presented in the following corollary.
\begin{cor}\label{cor_1}
    Consider system (\ref{sys_non_lin_def}) with a flat output $y=\varphi(x,u,u_{[1]}, \ldots)$. An $m$-tuple of time derivatives $y_{[\kappa]}=\varphi_{[\kappa]}(x,u,u_{[1]}, \ldots)$ of this flat output, with $\kappa \leq R$ and $\# \kappa = n$, can be introduced as new input $w$ by a quasi-static feedback law $u=\alpha(x,w,w_{[1]},\ldots)$ of the classical state if and only if the differentials $\D x, \D \varphi_{[\kappa]}, \D \varphi_{[\kappa+1]}, \ldots$ are linearly independent.
\end{cor}

The linear independence of $\D x, \D \varphi_{[\kappa]}, \D \varphi_{[\kappa+1]}, \ldots$ ensures that the trajectory of the newly chosen input $y_{[\kappa]}(t) = w(t)$ is not restricted by the current state $x$ in any form.
\vspace{-0pt}
\subsection{ Lagrangian control systems and configuration flatness }
\vspace{-12pt}
Considering Lagrangian systems, the $p$-dimensional smooth manifold $\Qm$ denotes the configuration space with local coordinates $(q^1, \ldots, q^p)$, and $\T(\Qm)$ denotes the tangent bundle over $\Qm$ with local coordinates $(q^1, \ldots, q^p$, $\dot{q}^1, \ldots, \dot{q}^p)$. When external control forces $u=(u^1,\ldots,u^m)$ are additionally applied, we speak of Lagrangian control systems. As, e.g., in \cite{SatoIwai:2012}, we only consider forces that act on the system through multiplication by a matrix $G(q)$. Hence, the equations of motion are given by
\begin{equation}
	\label{equ_of_motion}
	\frac{\D}{\D t} \left( \frac{\partial L(q,\dot{q})}{\partial \dot{q}^i} \right) - \frac{\partial L(q,\dot{q})}{\partial q^i} = G_{ij}(q)u^j,
\end{equation}
with $i=1,\ldots,p$ and $j=1,\ldots,m$. By solving (\ref{equ_of_motion}) for $\ddot{q}$, we obtain the classical state-space representation
\begin{equation}
	\label{lag_sys_class}
	\begin{aligned}
		\dot{q}^i & = v^i,  \\
		\dot{v}^i & = f^i(q,v,u), \quad i=1,\ldots,p,
	\end{aligned}
\end{equation}
with the $2p$-dimensional classical state $(q,v)$ and the $m$-dimensional input $u$. Based on the findings of \cite{RathinamMurray:1998}, we briefly review configuration flat Lagrangian control systems with $m=p-1$ inputs and a flat output of the form
\begin{equation}
	\label{flat_output}
	y^j = \varphi^j(q), \quad j=1,\ldots,p-1.
\end{equation}
The flat output (\ref{flat_output}) represents a smooth submersion \mbox{$\varphi: \Qm \rightarrow \mathcal{Y}$}, where $\mathcal{Y}$ denotes a smooth manifold with local coordinates $(y^1, \ldots, y^{p-1})$. Hence, we can formulate a local diffeomorphism of the form
\begin{equation}
	\label{trf_1}
	\begin{aligned}
		\bar{q}^j & = g^j(q) = \varphi^j(q), \quad j=1,\ldots,p-1, \\
		\bar{q}^p & = g^p(q),
	\end{aligned}
\end{equation}
and introduce transformed coordinates $\bar{q}$ on the configuration space $\Qm$. Based on (\ref{trf_1}), we obtain a regular state transformation of the form
\begin{subequations}
	\label{state_trf_1}
	\begin{align}
		\bar{q}^i & = g^i(q),  \label{state_trf_1_q} \\
		\bar{v}^i & = \partial_{q^k} g^i(q)v^k, \quad i,k=1,\ldots,p, \label{state_trf_1_v}
	\end{align}
\end{subequations}
on the tangent bundle $\T(\Qm)$. Hence, we can choose local coordinates $(y^1, \ldots, y^{p-1}, \bar{q}^p)$ on $\Qm$.
\vspace{-3pt}
\section{ Exact Linearization }\label{sec_5}
\vspace{-3pt}
In this section, we examine the exact linearization of configuration flat Lagrangian control systems (\ref{lag_sys_class}) with $m=p-1$ control inputs by a quasi-static feedback law of the classical state $(q,v)$ of the form
\begin{equation}
	\label{lin_fb_1}
	u = \alpha( q,v,w, w_{[1]}, \ldots )
\end{equation}
such that the closed-loop system 
\begin{equation}
	\label{closed_loop_sys}
	\begin{aligned}
		\dot{q} & = v, \\
		\dot{v} & = f\left(q,v,\alpha( q,v,w, w_{[1]}, \ldots )\right)\\
	\end{aligned}
\end{equation}
possesses a linear input-output behavior of the form (\ref{input_output}) between the new input $w$ and the flat output $y$. First, we present a detailed analysis of the flat parameterization of minimally underactuated configuration flat Lagrangian control systems. Subsequently, we demonstrate how to choose suitable orders $\kappa^j$ of the time derivatives of components $y^j$ of the flat output $y$ intended as new inputs. Based on the Jacobian matrix of the flat parameterization of the transformed coordinates $(\bar{q}, \bar{v})$, we demonstrate how to select suitable multi-indices $\kappa$.

\begin{thm} \label{thm_3}
    Consider a Lagrangian control system (\ref{lag_sys_class}) with a flat output $y=\varphi(q,v,u,u_{[1]},\ldots)$. The parameterization of the system variables $(q,v,u)$ by the flat output $y$ is of the form
	\begin{subequations}\vspace{-5pt}
		\label{para_map}
		\begin{align}
			\label{para_map_q} q & = F_q(y_{[0,R-2]}), \\
			\label{para_map_v} v & = F_v(y_{[0,R-1]}), \\
			\label{para_map_u} u & = F_u(y_{[0,R]}).
		\end{align}
	\end{subequations}
    For the special case with $m=p-1$ control inputs and a configuration flat output (\ref{flat_output}), the parameterizing map (\ref{para_map}) has the following properties:
    \begin{enumerate}[label=\roman*)]
        \item The multi-index $R=(r^1,\ldots,r^{p-1})$ meets \mbox{$2 \leq R \leq 4$}.
        \item The parameterization (\ref{para_map_q}) depends explicitly on the second-order time derivative of at least one \mbox{component} of the flat output (\ref{flat_output}).
        \item The third-order and fourth-order time derivatives of these components of the flat output (from property ii) appear explicitly in the parameterizations (\ref{para_map_v}) and (\ref{para_map_u}), respectively. 
    \end{enumerate}
\end{thm}
\vspace{-5pt}
\begin{pf}
    Let $R$ be the minimal multi-index as in Definition \ref{def_1}. Then, the parameterization of the system variables $(q,v,u)$ of a configuration flat Lagrangian control system (\ref{lag_sys_class}) is given by 
    \vspace{-5pt}
    \begin{equation} \label{para_map_def_0}
        \begin{aligned}
            q & = F_q(y_{[0,R-1]}), \\
            v & = F_v(y_{[0,R-1]}), \\
            u & = F_u(y_{[0,R]}). \\
        \end{aligned}
    \end{equation}
    Given the structure of (\ref{lag_sys_class}), we know that the parameterization $F_v$ of the generalized velocities is related to the parameterization $F_q$ of the generalized coordinates by
    \begin{equation} \label{calc_F_v_0}
        F^i_v = \sum^{r^j-1}_{\beta^j=0} \frac{\partial F^i_q }{\partial y^j_{[\beta^j]}} y^j_{[\beta^j+1]}, \vspace{-1pt}
    \end{equation}
    with $i=1,\ldots,p$ and $j=1,\ldots,p-1$. Since $R$ is the unique and minimal multi-index corresponding to the parameterizing map (\ref{para_map_def_0}), it follows that $F_v$ cannot explicitly depend on time derivatives of $y$ of the order $R$ or higher. Consequently, as indicated by (\ref{calc_F_v_0}), $F_q$ cannot explicitly depend on time derivatives of $y$ of the order $R-1$ or higher. However, $F_q$ must explicitly depend on $y_{[R-2]}$. Otherwise, $F_v$ would not depend on time derivatives of $y$ up to the order $R-1$, and thus $R$ would not be minimal. 

    We consider configuration flat Lagrangian control systems (\ref{lag_sys_class}) with $m=p-1$ inputs.
    Note that the multi-index $R$ is independent of the choice of coordinates on $\T(\Qm)$. Let us first substitute the parameterization (\ref{para_map_q}) of the generalized coordinates $q$ into the right-hand side of the diffeomorphism (\ref{trf_1}). Subsequently, we apply (\ref{calc_F_v_0}) for calculating the parameterization $F_{\bar{v}}$ of the transformed velocities based on $F_{\bar{q}}$. Accordingly, we obtain the flat parameterization of the transformed coordinates $(\bar{q}, \bar{v})$:
	\begin{equation}
		\label{flat_para_sec_2}
		\arraycolsep=3pt
		\begin{array}{rclcrcl}
			\bar{q}^1 &=& y^1, & & \bar{v}^1 &=& y^1_{[1]}, \\ [0.2cm]
			&\vdots&  & &  &\vdots&  \\ [0.2cm]
			\bar{q}^{p-1} &=& y^{p-1}, & & \bar{v}^{p-1} &=& y^{p-1}_{[1]}, \\ [0.2cm]
			\bar{q}^{p} &=& F^p_{\bar{q}} \left( y_{[0, R-2 ]} \right), & & \bar{v}^{p} &=& F^p_{\bar{v}} \left( y_{[0, R-1 ]} \right). \\ [0.2cm]
		\end{array}
            \vspace{-3pt}
	\end{equation}
    Next, let us prove the properties of (\ref{para_map}): \vspace{4pt} \break 
    i) First, we analyze the first column of (\ref{flat_para_sec_2}). The parameterization of the transformed coordinates $\bar{q}$ depends on all components of the flat output $y$. Given the diffeomorphism (\ref{trf_1}), the parameterization (\ref{para_map_q}) of the generalized coordinates $q$ must also depend explicitly on all components of $y$. Since the highest occurring time derivatives of $y$ in (\ref{para_map_q}) are of the order $R-2$, the lower bound for $R$ is given by $R\geq 2$. In \cite{KolarSchoberlSchlacher:2016-3-ifac-online} a bound for the maximum order of the time derivatives that occur in the parameterizing map of $(x, u, u_{[1]}, \ldots, u_{[\nu]})$-flat systems is presented. For the particular class of configuration flat systems, this bound is given by \mbox{$R \leq n - 2(m-1)$} with $n$ states and  $m$ control inputs. Hence, for configuration flat Lagrangian control systems (\ref{lag_sys_class}) with $n=2p$ states and $m=p-1$ control inputs, it follows that $R \leq 4$. \vspace{4pt} \break 
    ii) Because $(\bar{q}, \bar{v})$ are coordinates on $\T(\Qm)$, the differentials $\D \bar{q}$ and $\D \bar{v}$ must be linearly independent. Together with (\ref{flat_para_sec_2}), it follows that the flat parameterization $F^p_{\bar{q}}$ of the generalized coordinate $\bar{q}^p$ explicitly depends on the second-order time derivative of at least one component of the flat output $y$. Otherwise, $\D \bar{q}^p$ could be expressed by linear combination of the differentials $\D \bar{q}^1, \ldots, \D \bar{q}^{p-1}, \D \bar{v}^1, \ldots, \D \bar{v}^{p-1}$. Since $q$ and $\bar{q}$ are related by the diffeomorphism (\ref{trf_1}), (\ref{para_map_q}) explicitly depends on the second-order time derivative of at least one component of the flat output $y$. \vspace{4pt} \break 
    iii) From $R \leq 4$ and the fact that (\ref{para_map_q}) explicitly depends on the second-order time derivative of at least one component of $y$, it immediately follows that the highest occurring time derivatives of $y$ in (\ref{para_map_q}) are exactly of order two. Because of the structure of (\ref{lag_sys_class}), the third-order time derivatives of these components occur explicitly in (\ref{para_map_v}), and the fourth-order time derivatives of these components occur explicitly in (\ref{para_map_u}). \hfill $\Box$
\end{pf}

Before we demonstrate how to select suitable multi-indices $\kappa$, let us transfer Corollary \ref{cor_1} to Lagrangian control systems of the form (\ref{lag_sys_class}) and further simplify the condition that the differentials $\D q, \D v, \D \varphi_{[\kappa]}, \D \varphi_{[\kappa+1]}, \ldots$ must be linearly independent.
\begin{prop}\label{thm_4}
    Consider a Lagrangian control system (\ref{lag_sys_class}) with a flat output $y=\varphi(q,v,u,u_{[1]}, \ldots)$ as well as the corresponding parameterization (\ref{para_map}). An $m$-tuple of time derivatives $y_{[\kappa]}=\varphi_{[\kappa]}(q,v,u,u_{[1]}, \ldots)$ of this flat output with $\kappa \leq R$ and $\# \kappa = 2p$ can be introduced as new input $w$ by a quasi-static feedback law (\ref{lin_fb_1}) of the classical state $(q,v)$ if and only if the $2p\times 2p$ Jacobian matrix
    \begin{equation}
        \label{jacobian_0}
        \partial_{y_{[0,\kappa-1]}} 
        \begin{pmatrix}
            F_q \\ F_v
        \end{pmatrix}    \vspace{-3pt}
    \end{equation}
    is regular.
\end{prop}
\vspace{-3pt}
\begin{pf}
    According to Corollary \ref{cor_1}, the existence of a linearizing feedback law (\ref{lin_fb_1}) such that the closed loop system (\ref{closed_loop_sys}) possesses a linear input-output behavior of the form (\ref{input_output}) with $\# \kappa = 2p$ requires that the differentials
    \begin{equation}
    	\label{diff_1}
    	\D q, \D v, \D \varphi_{[\kappa]}, \D \varphi_{[\kappa+1]}, \ldots
    \end{equation}
    with time derivatives of $\varphi$ up to arbitrary order must be linearly independent. By utilizing the flat parameterization (\ref{para_map}), we express (\ref{diff_1}) in the form
    \begin{equation*}
        \label{diff_2}
        \begin{aligned}
            \D q^i & = \partial_{y^j} F^i_q \D y^j+\partial_{y^j_{[1]}}F^i_q \D y^j_{[1]} + \cdots  \\
            & \hspace{5mm} \cdots + \partial_{y^j_{[\kappa^j-1]}}F^i_q \D y^j_{[\kappa^j-1]} + \partial_{y^j_{[\kappa^j]}}F^i_q \D y^j_{[\kappa^j]} + \cdots, \\ 
            \D v^i & = \partial_{y^j} F^i_v \D y^j+\partial_{y^j_{[1]}}F^i_v \D y^j_{[1]} + \cdots  \\
            & \hspace{5mm} \cdots + \partial_{y^j_{[\kappa^j-1]}}F^i_v \D y^j_{[\kappa^j-1]} + \partial_{y^j_{[\kappa^j]}}F^i_v \D y^j_{[\kappa^j]} + \cdots, \\[-5pt] 
            \D \varphi^j_{[\kappa^j]} & = \D y^j_{[\kappa^j]},	
        \end{aligned}
    \end{equation*}
    with $i=1,\ldots,p$ and $j=1,\ldots,m$. By forming linear combinations, it follows that the linear independence of the differentials shown above is equivalent to the linear independence of the differential forms
    \begin{equation*}
        \label{diff_3}
        \begin{aligned}
            \partial_{y^j} F^i_q \D y^j + \partial_{y^j_{[1]}}F^i_q \D y^j_{[1]} + \cdots + \partial_{y^j_{[\kappa^j-1]}}F^i_q \D y^j_{[\kappa^j-1]}&, \\ 
            \partial_{y^j} F^i_v \D y^j + \partial_{y^j_{[1]}}F^i_v \D y^j_{[1]} + \cdots + \partial_{y^j_{[\kappa^j-1]}}F^i_v \D y^j_{[\kappa^j-1]}&,
        \end{aligned}
    \end{equation*}
    with $i=1,\ldots,p$ and $j=1,\ldots,m$. This further means that the rows of the Jacobian matrix of the parameterizations $F_q$ and $F_v$ with respect to time derivatives of the flat output $y$ up to the order $\kappa-1$ have to be linearly independent. Given $\# \kappa = 2p$, it follows that linear independence of (\ref{diff_1}) is equivalent to the regularity of (\ref{jacobian_0}). \hfill $\Box$
\end{pf}
We proceed by examining multi-indices $\kappa$ corresponding to feedback laws that allow for rest-to-rest transitions. Recall that we assume that (\ref{para_map}) is well defined around equilibrium points $y_s$. The following theorem presents all possible multi-indices $\kappa$ meeting the rank condition of Proposition~\ref{thm_4}, while excluding multi-indices that correspond to feedback laws that are guaranteed singular at equilibrium points.
\begin{thm}
    \label{thm_5}
    Considering a Lagrangian control system (\ref{lag_sys_class}) with $m=p-1$ control inputs and a configuration flat output (\ref{flat_output}) as well as the corresponding parameterization of the system variables (\ref{para_map}), the following applies:
    \begin{enumerate}[label=\roman*)] 
    \item For every multi-index $\kappa=(\kappa^1, \ldots, \kappa^{p-1})$ given by
    \begin{equation}
        \label{kappa_1}
        \kappa^j = 
        \begin{cases}
            4 \text{ for one } j \in \{ 1,\ldots,p-1 \} \text{ where } r^j = 4, \\
            2 \text{ for all other } j,
        \end{cases}
    \end{equation}
    there exists a quasi-static feedback law (\ref{lin_fb_1}) such that a linear input-output behavior (\ref{input_output}) of the closed-loop system (\ref{closed_loop_sys}) is achieved. Further, one can find at least one multi-index $\kappa$ of the form (\ref{kappa_1}).
    \item An input-output behavior (\ref{input_output}) with multi-indices $\kappa$ that are not of the form (\ref{kappa_1}) can either not be achieved by a quasi-static feedback law (\ref{lin_fb_1}) at all or only by a quasi-static feedback law (\ref{lin_fb_1}) which exhibits a singularity at equilibrium points.
    \end{enumerate}
\end{thm}
\begin{pf}
    According to Proposition \ref{thm_4}, a multi-index $\kappa$ is feasible if the submatrix (\ref{jacobian_0}) of the Jacobian matrix of $(F_q,F_v)$ is regular. Thus, we analyze which possibilities for constructing such regular submatrices exist. Since the linear independence of columns of this Jacobian matrix is not affected by state transformations, we can exploit the beneficial structure of the parameterizing map (\ref{flat_para_sec_2}) in transformed coordinates (\ref{state_trf_1}). In these coordinates, the Jacobian matrix of $(F_q, F_v)$ is given by
    \begin{equation}
        \setlength{\arraycolsep}{4pt}
        \label{jacobian_1}
        \begin{aligned}
        \partial_{y_{[0,3]}} 
        \begin{pmatrix}
            F_{\bar{q}} \\ F_{\bar{v}}
        \end{pmatrix} = 
        \left( \begin{array}{cccc}
            I & 0 & 0 & 0 \\
            \partial_y F^p_{\bar{q}} & \partial_{y_{[1]}} F^p_{\bar{q}} & \partial_{y_{[2]}} F^p_{\bar{q}} & 0  \\
            \cmidrule[0.4pt]{1-4}
            0 & I & 0 & 0 \\
            \undermat{A}{\partial_y F^p_{\bar{v}} & \partial_{y_{[1]}} F^p_{\bar{v}}} & \undermat{B}{\partial_{y_{[2]}} F^p_{\bar{v}} & \partial_{y_{[3]}} F^p_{\bar{v}}} \\
        \end{array} \right), 
        \end{aligned} \vspace{5pt}
    \end{equation} 
    with the $(p-1)\times(p-1)$-dimensional identity matrix $I$.\footnote{Note that we only have to consider the Jacobian matrix of $(F_q, F_v)$ with respect to consecutive time derivatives of $y$ up to order three, since higher time derivatives do not occur in the flat parameterizations (\ref{para_map_q}) and (\ref{para_map_v}).} \vspace{4pt} \break
    i) Since the $2p\times2(p-1)$ dimensional matrix $A$ has full rank $2(p-1)$, we choose two linearly independent columns of $B$ that can be appended to $A$ to obtain a regular submatrix (\ref{jacobian_0}). With the parameterization $(F_q, F_v)$ of the original coordinates and the relation (\ref{calc_F_v_0}), one can deduce that $\partial_{y_{[2]}} F_q = \partial_{y_{[3]}} F_v$. Similarly, we have $\partial_{y_{[2]}} F^p_{\bar{q}} = \partial_{y_{[3]}} F^p_{\bar{v}}$ for the parameterization $(F_{\bar{q}}, F_{\bar{v}})$ of the transformed coordinates. Given that fact, property ii) of Theorem \ref{thm_3} guarantees that \mbox{$\partial_{y^j_{[2]}} F^p_{\bar{q}} = \partial_{y^j_{[3]}} F^p_{\bar{v}} \neq 0$} applies with respect to at least one component $y^j$ of the flat output $y$. Thus, the rank of $B$ is exactly 2. Accordingly, appending two columns with the corresponding nonzero entries $\partial_{y^j_{[2]}} F^p_{\bar{q}}$ to $A$ yields a regular matrix (\ref{jacobian_0}) with a multi-index $\kappa$ of the form (\ref{kappa_1}). \vspace{4pt} \break
    ii) To obtain a feedback law (\ref{lin_fb_1}) that is well defined around equilibrium points, we have to choose multi-indices $\kappa$ such that the corresponding Jacobian (\ref{jacobian_0}) evaluated at $y_s$ is regular. Recall that the rank of $A$ is constant. By analyzing $B$ at equilibrium points $y_s$ it follows that
    \begin{equation} \label{equi_cons_2}
        \partial_{y_{[2]}} F^p_{\bar{v}} \Big|_{y_s} = 0.
    \end{equation}
    The proof for (\ref{equi_cons_2}) is given in Appendix \ref{appendix}. Therefore, the only possibility to construct a submatrix (\ref{jacobian_0}) that is regular at equilibrium points is to append two columns
    \begin{equation}\label{B_matrix}
    \begin{pmatrix}
        0 & 0 \\
        \partial_{y^j_{[2]}} F^p_{\bar{q}}|_{y_s} & 0 \\
        0 & 0 \\
        0 & \partial_{y^j_{[2]}} F^p_{\bar{q}}|_{y_s} \\
    \end{pmatrix}
    \end{equation}
    from the submatrix $B$ 
    to the submatrix $A$. Consequently, (\ref{jacobian_0}) can only be regular at equilibrium points if $\kappa$ is of the form (\ref{kappa_1}). \hfill $\Box$
\end{pf}
\begin{rem}
    Note, however, that multi-indices of the form (\ref{kappa_1}) can also possibly result in quasi-static feedback laws that are singular at equilibrium points.
\end{rem}
To derive a linearizing feedback law (\ref{lin_fb_1}), we proceed by substituting $y_{[\kappa, R-1]} = w_{[0,R-\kappa-1]}$ into the flat parameterizations (\ref{para_map_q}) and (\ref{para_map_v}), and obtain
\begin{equation}
	\begin{pmatrix}
		q \\ v
	\end{pmatrix} = 
	\begin{pmatrix}
		F_q(y_{[0,\kappa-1]}, w_{[0,R-\kappa-2]}) \\
		F_v(y_{[0,\kappa-1]}, w_{[0,R-\kappa-1]})
	\end{pmatrix}.
\end{equation}
Since the multi-index $\kappa$ is chosen such that (\ref{jacobian_0}) has full rank, the implicit function theorem applies. Hence, a solution of the form
\begin{equation}
	\label{y_kapp}
	y_{[0,\kappa-1]} = \psi (q, v, w_{[0,R-\kappa-1]})
\end{equation}
must exist locally. By substituting (\ref{y_kapp}) and $y_{[\kappa, R]} = w_{[0,R-\kappa]}$ into the flat parameterization (\ref{para_map_u}) of the input $u$, we finally obtain the linearizing feedback
\begin{equation}
	u = F_u(\psi (q, v, w_{[0,R-\kappa-1]}), w_{[0,R-\kappa]})
\end{equation}
of the desired form (\ref{lin_fb_1}).
\section{ Example }\label{sec_6}
In this section, we apply Theorem \ref{thm_5} to a Lagrangian control system, namely the well-known planar vertical take-off and landing (VTOL) aircraft with three degrees of freedom and two control inputs as discussed, e.g., in \cite{FliessLevineMartinRouchon:1999}, \cite{Gstottner:2022}, or \cite{Gstottner:2023-2}. The state-space representation (\ref{lag_sys_class}) is given by
\begin{equation*}
    \arraycolsep=3pt
    \begin{array}{rclcrcl}
        \dot{x} & = & v_x, & & \dot{v}_x & = & \epsilon\cos(\theta)u^2 - \sin(\theta)u^1, \\ [0.2cm]
        \dot{z} & = & v_z, & & \dot{v}_z & = & \cos(\theta)u^1 + \epsilon\sin(\theta)u^2-1, \\ [0.2cm]
        \dot{\theta} & = & \omega, & & \dot{\omega} & = & u^2,
    \end{array}
\end{equation*}
with the classical state $(q,v)=(x,z,\theta,v_x,v_z,\omega_\theta)$ and the control inputs $(u^1,u^2)$. The configuration flat output $y=(x-\epsilon\sin(\theta), z+\epsilon\cos(\theta))$ of the VTOL has been derived, e.g., in \cite{Martin:1996}, with (\ref{para_map}) given by
\begin{subequations}
	\label{para_vtol}
	\begin{align}
		\label{para_vtol_q_v} (q,v) & = \left(F_q(y^1, y^1_{[2]}, y^2, y^2_{[2]}), F_v(y^1_{[1,3]}, y^2_{[1,3]}) \right), \\
		\label{para_vtol_u} u & = F_u(y^1_{[2,4]}, y^2_{[2,4]}).
	\end{align}
\end{subequations} 
The multi-index denoting the highest occurring time derivatives in (\ref{para_vtol}) is $R=(4,4)$. A regular state transformation of the form (\ref{state_trf_1}) on $\T(\Qm)$ is given by
\begin{equation}\label{transf_coord}
	\arraycolsep=1.5pt
	\begin{array}{rclcrcl}
		\bar{q}^1 & =& x-\epsilon\sin(\theta), & \quad & \bar{v}^1 & = & v_x - \epsilon\omega\cos(\theta), \\
		\bar{q}^2 & =& z+\epsilon\cos(\theta), & \quad & \bar{v}^2 & = & v_z -\epsilon\omega\sin(\theta), \\
		\bar{q}^3 & =& \theta, & \quad & \bar{v}^3 & = & \omega. \\	
	\end{array}
\end{equation}
The flat parameterization $(F_{\bar{q}}, F_{\bar{v}})$ of the transformed coordinates (\ref{transf_coord}) is of the form (\ref{flat_para_sec_2}). According to item~i) of Theorem~\ref{thm_5}, the VTOL can be exactly linearized by a feedback law (\ref{lin_fb_1}) such that the closed loop system (\ref{closed_loop_sys}) shows a linear input-output behavior (\ref{input_output}) with the lengths of the decoupled integrator chains given by $\kappa_1=(4,2)$ and $\kappa_2=(2,4)$.
However, since we want to derive a feedback law (\ref{lin_fb_1}) that allows rest-to-rest transitions, we evaluate the Jacobian matrix (\ref{jacobian_1}) at equilibrium points and obtain
\begin{equation}
	\label{jacobian_vtol}
        \partial_{y_{[0,3]}} 
        \begin{pmatrix}
            F_{\bar{q}} \\ F_{\bar{v}}
        \end{pmatrix} \Big|_{y_s} = 
	\begin{pmatrix}
		1 & 0 & 0 & 0 & 0 & 0 & 0 & 0 \\
		0 & 1 & 0 & 0 & 0 & 0 & 0 & 0 \\
		0 & 0 & 0 & 0 & -1 & 0 & 0 & 0 \\
		0 & 0 & 1 & 0 & 0 & 0 & 0 & 0 \\
            0 & 0 & 0 & 1 & 0 & 0 & 0 & 0 \\
            0 & 0 & 0 & 0 & 0 & 0 & -1 & 0 \\
	\end{pmatrix}.
\end{equation}
Given (\ref{jacobian_vtol}), the only multi-index corresponding to a Jacobian (\ref{jacobian_0}) that is regular around $y_s$ is $\kappa_1=(4,2)$. Next, substituting $y_{[\kappa_1,R-1]}=(y^2_{[2,3]})=(w^2_{[0,1]})$ into (\ref{para_vtol_q_v}) yields
\begin{equation}\label{eq_fb_vtol_1}
	(q,v) = \left(F_q(y^1, y^1_{[2]}, y^2, w^2), F_v(y^1_{[1,3]}, y^2_{[1]}, w^2_{[0,1]})\right).
\end{equation}
According to the implicit function theorem, we can solve (\ref{eq_fb_vtol_1}) locally for $y_{[0,\kappa_1-1]}=(y^1_{[0,3]}, y^2_{[0,1]})$. Substituting the solution $y_{[0,\kappa_1-1]} = \psi(q, v, w^2_{[0,1]})$ as well as $y_{[\kappa_1,R]}=(y^1_{[4]}, y^2_{[2,4]})=(w^1, w^2_{[0,2]})$ into (\ref{para_vtol_u}) yields 
\begin{equation}\label{fb_vtol}
    \begin{aligned}
        u^1 & = \alpha^1(\theta, \omega, w^2), \\
        u^2 & = \alpha^2(\theta, \omega, w^1, w^2_{[0,2]}),
    \end{aligned}
\end{equation}
which is of the desired form (\ref{lin_fb_1}) and enables rest-to-rest transitions. Hence, we systematically determined which lengths of integrator chains are suitable by analyzing (\ref{jacobian_vtol}). The same feedback law was also derived in \cite{Rudolph2021} using Brunovský states. However, applying the derivation of (\ref{fb_vtol}) as presented in \cite{Rudolph2021} to other minimally underactuated Lagrangian systems is not straightforward and, as pointed out in \cite{Rudolph2021}, does not generally yield a linearizing feedback of the form (\ref{lin_fb_1}). Our systematic derivation, on the other hand, directly applies to any minimally underactuated configuration flat Lagrangian control system, and by construction, a linearizing quasi-static feedback law of the classical state is guaranteed.

\vspace{-6pt}
\bibliography{bibliography}

\appendix
\vspace{-12pt}
\section{} \label{appendix}
\vspace{-12pt} 
The proof that (\ref{equi_cons_2}) applies is based on the findings of \cite{RathinamMurray:1998}. Specifically, we focus on equation (3.16) of \cite{RathinamMurray:1998}, which was derived by using Riemannian geometry and by means of the coordinate transformation (\ref{trf_1}), i.e., local coordinates $(y^1, \ldots, y^{p-1}, \bar{q}^p)$ defined on $\Qm$. The authors explain that it is necessary for configuration flatness that every coefficient of $\dot{\bar{q}}^p$ and $\ddot{\bar{q}}^p$ in (3.16) is equal to zero. Using the notation established in Section 2, we represent equation (3.16) by  
\begin{equation}\label{equ_rath}
    \alpha_k(y, \bar{q}^p)y^k_{[2]} + \beta_{lk}(y, \bar{q}^p)y^l_{[1]}y^k_{[1]} + \gamma(y, \bar{q}^p) = 0,
\end{equation}
with $l,k=1,\ldots,p-1$. Solving (\ref{equ_rath}) for $\bar{q}^p$ yields a solution of the form $\bar{q}^p=F_{\bar{q}}^p(y_{[0,2]}) = \Tilde{F}_{\bar{q}}^p(y, y^l_{[1]}y^k_{[1]}, y_{[2]})$ that depends quadratically on first-order time derivatives of $y$. Since $F^p_{\bar{v}}$ is related to $\Tilde{F}_{\bar{q}}^p$ by (\ref{calc_F_v_0}), it follows that
\begin{equation*}\label{A2}
    \begin{aligned}
        \partial_{y^i_{[2]}} F_{\bar{v}}^p \Big|_{y_s}\hspace{-2mm}=\partial_{y^i_{[2]}}\hspace{-1mm}\left( \partial_{y^j}\Tilde{F}^p_{\bar{q}}y^j_{[1]} + \hspace{-1mm} \partial_{y^j_{[1]}}\Tilde{F}^p_{\bar{q}}y^j_{[2]} + \hspace{-1mm} \partial_{y^j_{[2]}}\Tilde{F}^p_{\bar{q}}y^j_{[3]} \right)\hspace{-1mm} \Big|_{y_s}\hspace{-2mm} &=\\
        \partial_{y^i_{[1]}} \hspace{-1mm} \left(\Tilde{F}^p_{\bar{q}}(y, y^l_{[1]}y^k_{[1]}, y_{[2]})\right) \hspace{-1mm} \Big|_{y_s} \hspace{-2mm} & =0,
    \end{aligned}
\end{equation*} 
with $i,j=1,\ldots,p-1$. Thus, the derivatives of $F^p_{\bar{v}}$ with respect to second-order time derivatives $y_{[2]}$ of the flat output vanish at equilibrium points.

\end{document}